\documentclass[english]{article}
\usepackage[T1]{fontenc}
\usepackage[latin1]{inputenc}
\usepackage{geometry}
\makeatletter

\usepackage{babel}
\makeatother
\begin{document}
\title{Cosine of angle  and center of mass of an operator}

\author{\normalsize Kallol Paul and Gopal Das}
\date{~}

\maketitle

\begin{abstract}
\noindent We consider the notion of real center of mass and total center of mass of a bounded linear operator relative to another bounded linear operator and explore their relation with cosine and total cosine of a bounded linear operator acting on a complex Hilbert space. We  give another proof of the Min-max equality and then generalize it using the notion of orthogonality of  bounded linear operators. We also illustrate with examples an alternative method of calculating the antieigenvalues and total antieigenvalues for finite dimensional operators.
\end{abstract}

\section{Introduction:}
 Suppose T is a bounded linear operator on a complex Hilbert space H with inner product ( , ) and norm $ \| ~\|$. A bounded linear operator T is called strongly accretive if $ Re \langle Tx,x \rangle \geq m_T > 0 $ for all unit vectors x in H. For a strongly accretive, bounded operator T on a Hilbert space \\
\[ \sup_{\|x\| = 1}~~ \inf _{\epsilon > 0} \| (\epsilon T - I ) x \|^2 =  \inf _{\epsilon > 0} ~~\sup_{\|x\| = 1}\| (\epsilon T - I ) x \|^2 ~.\]
This is the Min-max equality in operator trigonometry and was obtained by Gustafson\cite{5,9} in 1968, Asplund and Pt$\acute{a}$k\cite{1} in 1971.\\
The angle of an operator was introduced in 1968 by Gustafson\cite{6} while studying the problems
in the perturbation theory of semi-group generators in \cite{7,8}. The cosine
of the angle of $T$ was defined by Gustafson as follows:\begin{eqnarray}
cos\phi(T) & = & \inf_{\left\Vert Tf\right\Vert \neq0}\frac{Re(Tf,f)}{\left\Vert Tf\right\Vert \,\left\Vert f\right\Vert }.\label{eq:1}\end{eqnarray}
The properties of $\cos\phi\left(T\right)$ are dependent on the real
part of numerical range $W(T)$ of $T$.\\
\noindent The quantity $\cos\phi\left(T\right)$ has another interpretation
as the first antieigenvalue of $T$,\begin{eqnarray}
\mu_{1}\left(T\right) & = & \inf_{\left\Vert Tf\right\Vert \neq 0}\frac{Re(Tf,f)}{\left\Vert Tf\right\Vert \,\left\Vert f\right\Vert }.\label{eq:2}\end{eqnarray}
This concept was also introduced by Gustafson\cite{9} and studied by Gustafson\cite{10,11,12}, Seddighin\cite{19}, Gustafson and Seddighin\cite{15} and Gustafson and Rao\cite{13,14}. In \cite{2} we  studied the structure of  of antieigenvectors of a strictly accretive operator and in \cite{16} we studied total antieigenvalues of a bounded normal operator. In \cite{17} we introduced the notion of symmetric antieigenvalues. The notion of the cosine of angle of an operator has a connection with the min-max  equality according to Gustafson's theory[\cite{5} - \cite{12}] as
\[ sin T = \sqrt{ 1 - cos^2T} = \min_{\epsilon > 0 } \| \epsilon T - I \| ~.\]
While studying the norm of the inner derivation Stampfli\cite{20} proved that for
any bounded linear operator T there exists a complex scalar $ z_{0} $  such that
$ \| T- z_{0}I\| \leq \|T-zI\| $ for all complex scalar z. He defined $z_0 $ as the center
of mass ( or center) of T . In the Banach space B(H,H) for any two operators T and A , T is  orthogonal to A in the sense of James[4] iff $\|T+\lambda A \| \geq \|T\| $ for all scalars $ \lambda $. Thus if $z_0$ is center of mass of an operator T, then $ T - z_{0}I $ is orthogonal to I. We  studied in \cite{18} the the notion of orthogonality of two bounded linear operators T and A in B(H,H) and proved that T is orthogonal to A iff  $ \exists $ $ \{x_{n}\},  \|x_{n}\| = 1 $ such that $ (A^{*}Tx_{n},x_{n})  \rightarrow  0 $ and $ \|Tx_{n}\| \rightarrow  \|T\| $.We introduce the notion of real center of mass and total center of mass of an operator relative to another operator and explore their relation with cosine and total cosine of an operator. We also give an easy proof of the Min-max equality.\\ \\
\section{Center of mass of an operator :}
Before describing the definition of center of mass of an operator we first prove the following two results:\\
\noindent \textbf{Lemma 2.1} \label{L: lm1}\\
    \noindent  For any $ T,~A~\in~B(H)$ if $ |\epsilon | > 2 \frac{\|T\|}{\|A\|} ,~ \|A\| \ne 0,$
 then $ \| T \| < \|T -\epsilon A \| $.\\
\noindent \emph{\bf Proof.}
  \[ T,~A~\in~B(H), |\epsilon | > 2 \frac{\|T\|}{\|A\|} \Rightarrow \parallel T~-~\epsilon A \parallel \geq \mid \epsilon \mid~\parallel A \parallel~-~\parallel T \parallel 
 > 2 \frac{\|T\|}{\|A\|} \parallel A \parallel - \parallel T \parallel = \parallel T \parallel. \]
 Hence the proof. \\
\noindent \textbf{Lemma 2.2}  \label{L:lm2}\\
  \noindent  Let T and A be two bounded linear operators on a complex Hilbert space H. Then there exists
   a real scalar  $ \epsilon_{0} $ such that
   \[ \|T-\epsilon_{0}A\| \leq \|T-\epsilon A \| ~~~~~\forall \epsilon\in R. \]
\noindent \emph{\bf Proof.}
  We consider the function $ f:~G ~\rightarrow ~R $   defined by
     \[ f(\epsilon) = \|T -\epsilon A \| ~\forall \epsilon \in G \]
  where R is the set of real numbers and $ G \subseteq R $  is a closed
  interval centered at origin with radius $ 2 \frac{\|T\|}{\|A\|} $.\\
  Then it is clear that f is continuous, bounded and G is compact. So f attains its infimum.
  Thus there exists a $ \epsilon_{0} \in G$ such that  $f(\epsilon) \geq f( \epsilon_{0} )~~\forall
  \epsilon \in G \subseteq R $ which implies  $ \|T-\epsilon_{0}A\| \leq \|T-\epsilon A \|~~
  \forall \epsilon \in G .$ \\
 \noindent  If $ \epsilon \notin G $ then $ |\epsilon | > 2 \frac{\|T\|}{\|A\|} $
  and so $ \|T-\epsilon A \|~>~\parallel T \parallel$ by the   above lemma.
  Also $ f(0) \geq f( \epsilon_{0} ) $ i.e., $ \parallel T \parallel~\geq~\|T-\epsilon_{0}A\|.
  $ Hence  $\|T-\epsilon_{0}A\| \leq \|T-\epsilon A \| $ if $ \epsilon \notin G. $\\
  Thus $ \|T-\epsilon_{0}A\| \leq \|T-\epsilon A \|$ for all
    $\epsilon~ \in R $. Hence the proof.\\
\noindent Thus for any two bounded linear operators T and A there exists a real scalar $ \epsilon_0 $ such that
\[ \|T-\epsilon_{0}A\| \leq \|T-\epsilon A \|~~ \forall     ~\epsilon~ \in R .\]
We define $ \epsilon_0 $ as the \textbf{real center of mass} of  operator T relative to the operator A.\\
\noindent Likewise one can show that for any two bounded linear operators T and A there exists a complex scalar $ \lambda_0 $ such that
\[ \|T-\lambda_{0}A\| \leq \|T-\lambda A \|~~ \forall     ~\lambda~ \in C .\]
We define $ \lambda_0 $ as the \textbf{total center of mass} of  operator T relative to the operator A.  The real center of mass and total center of mass of operator T relative to the operator A is not always uniquely defined. For A = I, $ \lambda_0$ is the center of mass of T introduced by Stampfli in \cite{20}.
\section{Orthogonality of two operators T and A in B(H,H) :}
In the Banach space B(H,H) for any two operators T and A , T is  orthogonal to A in the sense of James\cite{4} iff $\|T+\lambda A \| \geq \|T\| $ for all scalars $ \lambda $ in C. We say that T is real orthogonal to A iff
$\|T+\epsilon A \| \geq \|T\| $ for all scalars $ \epsilon $ in R. Thus if $ \epsilon_0$ and $ \lambda_0 $ are real center of mass and total center of mass of T relative to A respectively then $ T - \epsilon_0 A $ is real orthogonal to A and $ T - \lambda_0 A $ is orthogonal to A respectively. We next prove two theorems which is going to characterize the real center of mass.\\
\noindent \textbf{Theorem 3.1} \label{T:th1}\\
\noindent The set $ W_{0}(A) =\{ \epsilon \in R ~:~ \exists ~\{x_{n}\}~\subset 
H,~ \|x_{n}\|=1,~ \| Tx_{n} \| \rightarrow \| T \| $  and $ Re( Tx_{n},Ax_{n} ) \rightarrow \epsilon $
\}
is non-empty, closed and convex.\\
\noindent \emph{\bf Proof.}
 Clearly $ W_{0}(A) $ is a non-empty subset.\\
 We now show that  $ W_{0}(A) $ is closed. Let $\epsilon_{n} \in W_{0}(A) $ and $
   \epsilon_{n}~\rightarrow~\epsilon $. As $\epsilon_{n}\in W_{0}(A)$
   so there exists $ \{x_{k}^{n}\}_{k=1}^{\infty}
   , \parallel x_{k}^{n}\parallel ~=~1$
   for all k and n such that for each
   $n~=~1,2,3,......$
   \[ Re(A^{*}Tx_{k}^{n},x_{k}^{n}) ~\longrightarrow~\epsilon_{n}
   ~as~
   k~\longrightarrow~\infty \]
   and \[\|Tx_{k}^{n}\|~\longrightarrow~\|T\|~~as ~ k~\longrightarrow~\infty. \]
   Now for each n there exists $ k_{n}$ such that\\
   $\mid Re(A^{*}Tx_{k}^{n},x_{k}^{n})
   ~-~\epsilon_{n}\mid~<~\frac{1}{n}$ and
   $\mid \|Tx_{k}^{n}\|~-~\|T\| \mid~<~\frac{1}{n}$ for all $
   k~\geq~ k_{n}.$
   So we get \\
   $\mid Re(A^{*}Tx_{k_{n}}^{n},x_{k_{n}}^{n})
   ~-~\epsilon \mid~\leq~\mid Re(A^{*}Tx_{k_{n}}^{n},x_{k_{n}}^{n})
   ~-~\epsilon_{n}\mid~+~\mid \epsilon_{n} - \epsilon \mid
   ~\longrightarrow~0 $ as $n~\longrightarrow~\infty$.\\
   So $ Re(A^{*}Tx_{k_{n}}^{n},x_{k_{n}}^{n})~\longrightarrow~\epsilon $ as
   $n~\longrightarrow~\infty$.\\
   Similarly $ \|Tx_{k_{n}}^{n}\|~\longrightarrow~\|T\| $ as $n~\longrightarrow~\infty.$
    Hence $\epsilon  \in W_{0}(A). $ Thus $ W_{0}(A)$ is closed.\\
    Following an idea of Das et al\cite{3} we next show that $W_{0}(A)$ is convex. Let $ \mu $ and
    $  \eta \in  W_{0}(A)$ , $ \mu \neq \eta $ and t be any real scalar
   $ 0<t<1 $.  Then there exist $  \{x_{n}\}$ and $\{y_{n}\} $  in H with  $\|x_{n}\| = 1, \|y_{n}\| =
   1~ $ for all n , such that \[ Re(A^{*}Tx_{n},x_{n}) \rightarrow \mu ~ and ~  \|Tx_{n}\| \rightarrow
    \|T\|,\]
  \[ Re(A^{*}Ty_{n},y_{n}) \rightarrow \eta ~ and ~  \|Ty_{n}\| \rightarrow \|T\|.\]
For any real scalar $\epsilon$ and for any subsequence
$\{x_{n_{k}}\},~\{y_{n_{k}}\}$ of sequences $
\{x_{n}\},~\{y_{n}\}$ respectively. Suppose
\[ \lim_{k \longrightarrow \infty} \parallel x_{n_{k}}~\pm~\epsilon y_{n_{k}}\parallel ~=~ 0  \]
$\Rightarrow~ \lim_{k \longrightarrow \infty} \{
\parallel x_{n_{k}} \parallel^{2}~+~\epsilon^{2}\parallel y_{n_{k}}
\parallel^{2}~\pm~2\epsilon Re(x_{n_{k}},y_{n_{k}})\} ~=~ 0.$\\
$\Rightarrow~\lim_{k \longrightarrow \infty} Re(x_{n_{k}},y_{n_{k}}) ~=~ \mp~\frac{1+\epsilon^{2}}{2\epsilon} $\\
$\Rightarrow~\epsilon^{2}~ =~ 1$, otherwise $ \mid {\frac{1+ \epsilon^2}{2 \epsilon}} \mid > 1 $ which contradicts the fact that 
 $ \mid Re(x_{n_{k}},y_{n_{k}}) \mid  \leq 1 $.\\
Let $\delta_{n_{k}}~=~x_{n_{k}}~\pm~\epsilon y_{n_{k}}.$ Then
$\delta_{n_{k}}~\longrightarrow~0$ as $k~\longrightarrow~\infty.$\\
Now
\[\begin{array}{l}
 \mu~=~\lim_{k \longrightarrow
\infty}Re(A^{*}Tx_{n_{k}},x_{n_{k}}) \\
~~~=~\lim_{k \longrightarrow
\infty}Re(A^{*}T(\delta_{n_{k}}~\mp~\epsilon
y_{n_{k}}),\delta_{n_{k}}~\mp~\epsilon y_{n_{k}}) \\
~~~=~\epsilon^{2}\lim_{k \longrightarrow
\infty}Re(A^{*}Ty_{n_{k}},y_{n_{k}})\\
~~~=~\epsilon^{2} \eta \\
~~~ =~\eta. \\
\end{array} \]
 Thus we can assume that for any real scalar $\epsilon$ and
for any subsequence $\{x_{n_{k}}\},~\{y_{n_{k}}\}$ of sequences $
\{x_{n}\},~\{y_{n}\}$ respectively
\[ \lim_{k
\longrightarrow \infty}
\parallel x_{n_{k}}~\pm~\epsilon y_{n_{k}}\parallel~\neq~0. \]
Now, $\{Re(A^{*}T(x_{n}+\epsilon y_{n}),x_{n}+\epsilon y_{n})\}$ is a
bounded sequence and so it has a convergent subsequence, say,
$\{Re(A^{*}T(x_{n_{k}}+\epsilon y_{n_{k}}),x_{n_{k}}+\epsilon
y_{n_{k}})\}$. Also $ \{\parallel x_{n_{k}}~+~\epsilon
y_{n_{k}}\parallel \}, $ being a bounded sequence, has a
convergent subsequence, say $ \{\parallel x_{n_{k}}^{'}~+~\epsilon
y_{n_{k}}^{'}\parallel \} $. We show that there exists a real
scalar $\epsilon$ for which\\
$$
   \lim_{k \rightarrow \infty} \frac{Re(A^{*}T(x_{n_{k}}^{'} + \epsilon y_{n_{k}}^{'}),
  x_{n_{k}}^{'} + \epsilon y_{n_{k}}^{'})}{ \|x_{n_{k}}^{'} + \epsilon y_{n_{k}}^{'}\|^{2}}
   = t\mu + (1-t)\eta
 $$ 
\[\begin{array}{l}
  \Leftrightarrow \lim_{k \rightarrow \infty} [
  Re(A^{*}Tx_{n_{k}}^{'},x_{n_{k}}^{'})~+~\ \epsilon^{2}
  Re(A^{*}Ty_{n_{k}}^{'},y_{n_{k}}^{'})~+~ \epsilon \{
  Re(A^{*}Tx_{n_{k}}^{'},y_{n_{k}}^{'}) +  Re(A^{*}Ty_{n_{k}}^{'},x_{n_{k}}^{'}) \}\\\\
   ~-~\{t\mu + (1-t)\eta\}\|x_{n_{k}}^{'} + \epsilon
  y_{n_{k}}^{'}\|^{2}]~=~0\\\\
   \Leftrightarrow \mu~+~\epsilon^{2} \eta~+~\epsilon  \lim_{k \rightarrow
  \infty}  \{Re(A^{*}Tx_{n_{k}}^{'},y_{n_{k}}^{'}) + Re(A^{*}Ty_{n_{k}}^{'},x_{n_{k}}^{'})\}\\\\
  ~=~\{t\mu + (1-t)\eta\} \{1 +  \epsilon^{2}+~2\lim_{k \rightarrow
  \infty}Re \epsilon   (x_{n_{k}}^{'},y_{n_{k}}^{'}) \}
  \end{array}\]
This  on simplification yields
\[
\epsilon^{2}~+~C\epsilon~~-~\frac{1-t}{t}~=~0,~~~~~~~~~~(3)\]
where C is a real constant independent of $\epsilon $.\\
Equation(3) has  two nonzero roots of different signs, say 
$r_{1}$ and $-r_{2}$ ( $ r_1,~r_2 > 0 )$.\\
Now $ \{\parallel T(x_{n_{k}}^{'}~\pm~r
y_{n_{k}}^{'})\parallel \} $ is a bounded sequence and so it has a
convergent subsequence, say, $ \{\parallel T(
\tilde{x}_{n_{k}}~\pm~r \tilde{y}_{n_{k}})\parallel \}
$ where r is any positive real
number.\\
So,
\[\begin{array}{l}
\|T\|^2 \geq \lim_{k \rightarrow \infty} \frac{ \parallel T
(\tilde{x}_{n_{k}}~\pm~r \tilde{y}_{n_{k}})\parallel
^{2}}{ \parallel \tilde{x}_{n_{k}}~\pm~r
\tilde{y}_{n_{k}}\parallel^{2}}\\
~~~~~=~\parallel T \parallel ^{2}~\pm~2r~\lim_{k \rightarrow
\infty}~\frac {Re(T^{*}T\tilde{x}_{n_{k}}~-~\parallel T
\parallel
^{2}\tilde{x}_{n_{k}},\tilde{y}_{n_{k}})}{\parallel
\tilde{x}_{n_{k}}~\pm~r\tilde{y}_{n_{k}}\parallel^{2}}.
\end{array}\]
The sign of the second term on the R.H.S. is independent of r and
so for all positive real r,\\
either
\[ \parallel T (\frac{\tilde{x}_{n_{k}}~+~r
\tilde{y}_{n_{k}}}{\parallel \tilde{x}_{n_{k}}~+~r
\tilde{y}_{n_{k}}\parallel} )\parallel
~\longrightarrow~\parallel T \parallel  ~as
~k~\longrightarrow~\infty \] or
\[ \parallel T (\frac{\tilde{x}_{n_{k}}~-~r
\tilde{y}_{n_{k}}}{\parallel \tilde{x}_{n_{k}}~-~r
\tilde{y}_{n_{k}}\parallel} ) \parallel
~\longrightarrow~\parallel T \parallel  ~as
~k~\longrightarrow~\infty. \]
Let
\[ z_{n_{k}}~=~\frac{\tilde{x}_{n_{k}}~+~r
\tilde{y}_{n_{k}}}{\parallel \tilde{x}_{n_{k}}~+~r
\tilde{y}_{n_{k}}\parallel}  \]
and
\[ w_{n_{k}}~=~ \frac{\tilde{x}_{n_{k}}~-~r
\tilde{y}_{n_{k}}}{\parallel \tilde{x}_{n_{k}}~-~r
\tilde{y}_{n_{k}}\parallel}. \]
From above we see that either  $ \|Tz_{n_{k}}\|~\longrightarrow~\|T\| $ and $
Re(A^{*}Tz_{n_{k}},z_{n_{k}})~\longrightarrow~t\mu +
(1-t)\eta~as~k~\longrightarrow~\infty,$ or
$\|Tw_{n_{k}}\|~\longrightarrow~\|T\| $
 and $Re(A^{*}Tw_{n_{k}},w_{n_{k}})~\longrightarrow~t\mu +
(1-t)\eta~as~k~\longrightarrow~\infty.$\\
Hence $ t\mu + (1-t)\eta~\in~W_{0}(A)$. Thus $ W_{0}(A)$ is
convex. \\
This completes the proof of the theorem.\\
\noindent We now prove the following Theorem :\\
\noindent \textbf{Theorem 3.2} \label{T:th2}\\
 \noindent  $ \|T \| \leq  \|T - \epsilon A \| $
$ \forall \epsilon  \in $ R  iff $ \exists $ $ \{x_{n}\},  \|x_{n}\| = 1 $ such that $
Re(A^{*}Tx_{n},x_{n})  \rightarrow  0 $ and $ \|Tx_{n}\| \rightarrow
\|T\| $.\\
\noindent \emph{\bf Proof.}
 Suppose $ \exists~ \{x_{n}\}, \|x_{n}\| = 1 $ such that $ (A^{*}Tx_{n},x_{n})
\rightarrow 0  $ and $ \|Tx_{n}\| \rightarrow \|T\| $.\\
\noindent  Then $ {\|(T - \epsilon A)x_{n}\|}^{2} ={\|Tx_{n}\|}^{2} +{\epsilon}^{2}
{\|Ax_{n}\|}^{2} $ - 2Re $ \epsilon (Tx_{n},Ax_{n}) $, $ \forall \epsilon
\in $ R. So $ {\|T-\epsilon A \|}^{2} \geq \limsup_{n \rightarrow \infty}
{\|(T-\epsilon A)x_{n}\|}^{2} \geq {\|T\|}^{2} ~~ \forall \epsilon \in $ R.
Hence $ \|T-\epsilon A \| \geq \|T\|~~ \forall \epsilon \in $ R.  \\
\noindent Conversely let $ \|T\| \leq \|T- \epsilon A \| $ $ \forall \epsilon \in $ R. We need
to show $ 0 \in W_{0}(A) $.\\
\noindent Without loss of generality we can assume that $ \|A\| = 1 $.\\
Suppose $ 0 \not\in W_{0}(A) $ . Then as $  W_{0}(A) $ is closed and convex
by rotating T suitably we can assume that  $ W_{0}(A)  >  \eta > 0 $.\\
\noindent Let M = \{$x\in H ~:~  \|x\|$ = 1 and Re $ (Tx,Ax) \leq \eta $/2\} and
$\beta $ = $ sup_{x \in M} \|Tx\| $. We first claim that  $ \beta < \|T\| $.
Suppose $ \beta = \|T\|$. Then there exists $ x_n \in M $ such
that  $ \| Tx_n \| \rightarrow \|T\| $. As $ x_n \in M $ so $ Re (Tx_n,Ax_n)
\leq \eta/2 $ and $ \|x_n \|=1$. Now $ \{ Re(Tx_n,Ax_n)\} $ is a bounded sequence
 and so it has a convergent subsequence, without loss of generality we can assume
  that $ \{ Re(Tx_n,Ax_n)\} $ is convergent and converges to some point $ \mu $ (say).
   Then $ \mu \in W_0(A)$. Now $ Re(Tx_n,Ax_n) \leq \eta/2 $ and so $  \mu \leq \eta / 2$.
    This contradicts the fact that $  W_{0}(A) > \eta $.\\
\noindent Let $ \epsilon_0 = \min \{ \eta, \frac{\|T\| - \beta}{2 \|A\|} \} $. \\
 Let $ x \in M $ . Then $ \|(T-\epsilon_{0}A)x\| \leq \|Tx\| + \mid \epsilon_{0} \mid \|Ax\|  \leq \beta + \{( \|T\| - \beta )/ (2 \|A\|)  \} \|A\| = \|T\|/2  + \beta /2 $. So $ \sup_{x \in M ~and~ \|x\|=1 } \| (T-\epsilon_0 A)x \| \leq (\|T\| + \beta)/2 < \|T\| $. \\
Again let $ x \notin M $. Then let $ Tx = (a+ib)Ax + y $, where $ (Ax,y) $ = 0.
Now $ 2a \geq 2a ||Ax\|^2 = 2 Re (Tx,Ax) > \eta \geq \epsilon_0 $ and so $ 2a - \epsilon_0 > 0 $.
So
\begin{eqnarray*}
 {\|(T-\epsilon_{0}A)x\|}^{2}&=&\{{(a-\epsilon_{0})}^{2} + b^{2}\} {\|Ax\|}^{2}+{\|y\|}^{2}\hspace{2cm}      \\
                               &=&\mbox{}{\|Tx\|}^{2} + ({\epsilon_{0}}^{2} - 2a \epsilon_{0}) {\|Ax\|}^{2}\\
                                & \leq &  {\|Tx\|}^{2} + ({\epsilon_{0}}^{2} - 2a \epsilon_{0})  \\
                                & \leq &  {\|T\|}^{2} + ({\epsilon_{0}}^{2} - 2a \epsilon_{0})
\end{eqnarray*}	
Hence  	 $ \sup_{x \notin M ~and~ \|x\|=1 }  \| (T-\epsilon_0 A)x \|^2 \leq \|T\|^2 	+ ({\epsilon_{0}}^{2} - 2a \epsilon_{0}) $.	
Thus in all cases $ \| T - \epsilon_0 A \| < \|T\| $. Hence $ \exists \hspace{.15cm} \{x_{n}\}, \|x_{n}\|$ = 1 $ \forall n,~  (Tx_{n},Ax_{n}) \rightarrow 0 $ and $ \|Tx_{n}\| \rightarrow \|T\| $.\\
So far we proved that for any two operators T and A in B(H) with $ \| A\| \leq 1  $, $ \|T \| \leq  \|T - \epsilon A \| $ $ \forall \epsilon  \in $ R  implies that $ \exists $ $ \{x_{n}\},  \|x_{n}\| = 1 $ such that $
Re(A^{*}Tx_{n},x_{n})  \rightarrow  0 $ and $ \|Tx_{n}\| \rightarrow  \|T\| $.\\
\noindent   This completes the proof.\\
\noindent From the last theorem it follows that if $ \epsilon_0 $ is the real center of mass of T relative to A then
$ \exists $ $ \{x_{n}\},  \|x_{n}\| = 1 $ such that $ Re((T - \epsilon_0 A)x_{n},A x_{n})  \rightarrow  0 $ and $ \|(T - \epsilon_0 A)x_{n}\| \rightarrow  \|T - \epsilon_0 A\| $.\\
\noindent We next show the uniqueness of $ \epsilon_0 $ under the assumption that the approximate point spectrum of A, $\sigma_{app}(A)$ does not contain 0. Suppose
\[ \| T \| = \| T - \epsilon_0 A\| \leq \|T- \epsilon A \| ~ \forall ~\epsilon ~\in R ~~and~~ \epsilon_{0} \neq 0.\]
Then $ \exists \hspace{.15cm}\{x_{n}\},\|x_{n}\|$ = 1 such that
$ ( (T - \epsilon_{0}A)x_{n}, Ax_{n} ) \rightarrow 0 $ and
$ \|(T - \epsilon_{0}A)x_{n} \| \rightarrow \|T - \epsilon_{0}A \| $ . \\
So
\begin{eqnarray*}
{\|T - \epsilon_{0} A\|}^2 &=& \lim \{ {\|(T - \epsilon_{0} A)x_{n}\|}^{2} \} \\
                          &=& \lim \{ {\|Tx_{n}\|}^{2} +\epsilon_{0}^{2}{\|Ax_{n}\|}^{2}
                          - 2Re(\epsilon_{0} (Tx_{n},Ax_{n})) \}\\
                          &=& \lim \{ {\|Tx_{n}\|}^{2} -\epsilon_{0}^{2}{\|Ax_{n}\|}^{2} \} \\
                          &=& \lim \{ {\|Tx_{n}\|}^{2}\} -{\mid\epsilon_{0}\mid}^{2} \lim \{{\|Ax_{n}\|}^{2} \} \\
                          &<& \lim \{ {\|Tx_{n}\|}^{2}\}~~,~~~~~ since~0~~\notin~~ \sigma_{app}(A)\\
                          &\leq& { \|T\|}^2
\end{eqnarray*}
This contradicts the fact that $ \| T \| = \| T - \epsilon_0 A\| $. Hence $ \epsilon_0 = 0 $.  \\
We give an example to show that $ \epsilon_{0} $ may not be unique if $ 0~~\in~~ \sigma_{app}(A) $.\\
Let T and A be two bounded linear operators defined on $ R^{2} $ as  T (x,y) = (x,0) and A(x,y) = (0,y),  $ \forall ~(x,y) \in R^{2} $.  Then $ \| T \| = \| T - A \| = \| T - (-1) A \| \leq \| T  - \epsilon A \| ~~ \forall ~\epsilon \in R $.\\
\noindent We are now in a position to prove\\
\noindent \textbf{Theorem 3.3} \label{T:th3} \\
 \noindent \textbf{(Min-max equality)} For a strongly accretive, bounded operator T on a Hilbert space \\
\[ \sup_{\|x\| = 1}~~ \inf _{\epsilon > 0} \| (\epsilon T - I ) x \|^2 =  \inf _{\epsilon > 0} ~~\sup_{\|x\| = 1}\| (\epsilon T - I ) x \|^2 ~.\]
\noindent \emph{\bf Proof.}
For a fixed but arbitrary $ y \in H,~ \|y\| = 1 $ we have\\
\begin{eqnarray*}
\| (\epsilon T - I ) y \|^2 & \leq  &~~\sup_{\|x\| = 1}\| (\epsilon T - I ) x \|^2 \\
\Rightarrow \inf _{\epsilon > 0} \| (\epsilon T - I ) y \|^2 & \leq & \inf _{\epsilon > 0} ~~\sup_{\|x\| = 1}\| (\epsilon T - I ) x \|^2
\end{eqnarray*}
This holds for all y in H with $ \|y\| = 1 $. So
\[ \sup_{\|x\| = 1}~~ \inf _{\epsilon > 0} \| (\epsilon T - I ) x \|^2 \leq \inf _{\epsilon > 0} ~~\sup_{\|x\| = 1}\| (\epsilon T - I ) x \|^2 ~.\]
For the converse part we use the notion of center of mass of an operator. As T is strongly accretive   so there exists a real scalar $ \epsilon_0 > 0 $ such that
\[ \|  \epsilon_0 T - I \| = \inf _{\epsilon > 0} \| (\epsilon T - I ) \|.\]
Further there exists a sequence  $ \{x_{n}\},  \|x_{n}\| = 1 $ such that $ Re((\epsilon_0 T - I)x_{n}, Tx_{n})  \rightarrow  0 $ and $ \|(\epsilon_0 T - I)x_{n}\| \rightarrow  \|\epsilon_0 T - I\| $. \\
Now
\begin{eqnarray*}
\| \epsilon_0 T - I \|^2 & = & \lim_{n\rightarrow \infty} \|(\epsilon_0 T - I)x_n \|^2 \\
                         & = & \lim_{n\rightarrow \infty} \{ 1 - 2 \epsilon_0 Re \langle T x_n, x_n \rangle + {\epsilon_0}^2 \|Tx_n\|^2 \} \\
                         & = & \lim_{n\rightarrow \infty} \{ 1 - \frac{{ Re \langle T x_n, x_n \rangle }^2}{ \|Tx_n\|^2 } \}\\
                         & \leq & \sup_{\|x\| = 1} \{ 1 - \frac{{ Re \langle T x, x \rangle }^2}{ \|Tx\|^2 } \}~~~~~(*)
\end{eqnarray*}
For a fixed y in H with $ \|y \| = 1 $ we see that
\[ \| (\epsilon T - I ) y \|^2 = 1 - 2 \epsilon Re \langle T y, y \rangle + {\epsilon}^2 \|Ty\|^2 \]
achieves its minimum at $ \epsilon(y) = \frac{{ Re \langle T y, y \rangle }}{ \|Ty\|^2 } $ and the minimum value is
\[ \inf_{\epsilon > 0 } \| (\epsilon T - I ) y \|^2 = 1 - \frac{{ Re \langle T y, y \rangle }^2}{ \|Ty\|^2 } .\]
Using this in (*) we get
\[ \| \epsilon_0 T - I \|^2 \leq \sup_{\|x\| = 1} \inf_{\epsilon > 0 } \| (\epsilon T - I ) x\|^2 \]
and so
\[ \inf _{\epsilon > 0} ~~\sup_{\|x\| = 1}\| (\epsilon T - I ) x \|^2 \leq \sup_{\|x\| = 1}~~ \inf _{\epsilon > 0} \| (\epsilon T - I ) x \|^2 .\]
Thus we have the Min-max equality
\[ \sup_{\|x\| = 1}~~ \inf _{\epsilon > 0} \| (\epsilon T - I ) x \|^2 = \inf _{\epsilon > 0} ~~\sup_{\|x\| = 1}\| (\epsilon T - I ) x \|^2 ~.\]
\noindent Next, if $\epsilon_{0}$ is the real center of mass of a bounded linear operator $T$ , then there exists a sequence $\{ x_{n} \}$ in H, $\|x_{n}\|=1$, $Re((I-\epsilon_{0}T)x_{n},Tx_{n}) \rightarrow 0$ and $\|(I-\epsilon_{0}T)x_{n}\| \rightarrow \|I-\epsilon_{0}T\|$.\\
\noindent  For a bounded linear operator T the antieigenvalue is defined as \\
\[  \cos  T  =  \inf_{\| Tx \| \neq 0 } \frac{ Re \langle Tx,x \rangle } { \|Tx\| \|x\|}~.\]
\noindent We next prove the theorem\\
\noindent \textbf{Theorem 3.4} \label{T:th4} \\
\noindent Suppose $ \epsilon_{0} $ is the real center of mass of a strictly accretive, bounded linear operator T. Then $\cos T = \lim_{n \rightarrow \infty} \frac{Re(Tx_{n},x_{n})}{\|Tx_{n}\|}$ where $\{ x_{n} \}$ is a sequence of unit vectors in H, $Re((I-\epsilon_{0}T)x_{n},Tx_{n}) \rightarrow 0$ and $\|(I-\epsilon_{0}T)x_{n}\| \rightarrow \|I-\epsilon_{0}T\|$.\\
\noindent \emph{\bf Proof.}
 As $ \epsilon_0  $ is the real center of mass of T so there exists a sequence  $ \{x_{n}\},  \|x_{n}\| = 1 $ such that $ Re((I - \epsilon_0 T )x_{n}, Tx_{n})  \rightarrow  0 $ and $ \|(I - \epsilon_0 T )x_{n}\| \rightarrow  \|I - \epsilon_0 T\| $ and
\[ \|I - \epsilon_0 T\| = \inf_{\epsilon > 0} \| I - \epsilon T \| = \inf _{\epsilon > 0} ~~\sup_{\|x\| = 1}\| (I - \epsilon T  ) x \| ~.\]
As in Theorem 3.3  we have
\begin{eqnarray*}
\| \epsilon_0 T - I \|^2 & = & \lim_{n\rightarrow \infty} \{ 1 - \frac{{ Re \langle T x_n, x_n \rangle }^2}{ \|Tx_n\|^2 } \}\\
                         & \leq & \sup_{\|x\| = 1} \{ 1 - \frac{{ Re \langle T x, x \rangle }^2}{ \|Tx\|^2 } \} \\
                         & \leq & \sup_{\|x\| = 1} \inf_{\epsilon > 0 } \| (I - \epsilon T  ) x \|^2 \\
                         & = & \| \epsilon_0 T - I \|^2
\end{eqnarray*}
This shows that \[ \inf_{\|x\| = 1 }  \frac{{ Re \langle T x, x \rangle }}{ \|Tx\| } = \lim_{n\rightarrow \infty} \frac{{ Re \langle T x_n, x_n \rangle }}{ \|Tx_n\| } = \cos T.\]
\noindent We next give an example to calculate antieigenvalue for a finite dimensional operator using the above method.\\
\noindent \textbf{Example 3.5} \\
$T=\left(
  \begin{array}{cc}
    1 & 0 \\
    0 & 1+i \\
  \end{array}
\right)
$ be an operator on a two dimensional complex Hilbert space H. Let, $z=(z_{1}, z_{2})^{t} \in H$, where $\left|z_{1}\right|^{2} + \left| z_{2}\right|^{2} = 1$. Then, $Re(Tz,z)=1$ and $\|Tz\|= \sqrt{1+\left|z_{2}\right|^{2}}$. Now, $\sup_{\|z\| = 1}\inf_{\epsilon >0}\|(\epsilon T - I)z\| =\sqrt{0.5}$ and this supremum is attained by the vector $z_{0}=(z_{1},0)^{t}$, where $\left|z_{1}\right|=1$. Then, $\epsilon_{0}=0.5$ and $\|\epsilon_{0}T-I\|=\sqrt{0.5}$. Now, $\cos T = \frac{Re(Tz_{0},z_{0})}{\|Tz_{0}\|}=\frac{1}{\sqrt{2}}$.\\
\noindent In \cite{18} we proved that  if  $ \lambda_0  $ is the total center of mass of T then there exists a sequence  $ \{x_{n}\},  \|x_{n}\| = 1 $ such that $ ((I - \lambda_0 T )x_{n}, Tx_{n})  \rightarrow  0 $ and $ \|(I - \lambda_0 T )x_{n}\| \rightarrow  \|I - \lambda_0 T\| $.\\
\noindent  For a bounded linear operator T the total antieigenvalue is defined as \\
\[ \mid \cos \mid T  =  \inf_{\| Tx \| \neq 0 } \frac{ \mid \langle Tx,x \rangle \mid } { \|Tx\| \|x\|}~.\]
\noindent In \cite{16} we also studied the total antieigenvalue of a bounded linear operator. \\
\noindent We now prove the theorem \\
\noindent \textbf{Theorem 3.6} \\
 \noindent For a  bounded linear operator T on a Hilbert space \\
\[ \sup_{\|x\| = 1}~~ \inf _{\lambda \in C} \| (\lambda T - I ) x \|^2 =  \inf _{\lambda \in C} ~~\sup_{\|x\| = 1}\| (\lambda T - I ) x \|^2 ~.\]
\noindent \emph{\bf Proof.}
For a fixed but arbitrary $ y \in H,~ \|y\| = 1 $ we have\\
\begin{eqnarray*}
\| (\lambda T - I ) y \|^2 & \leq  &~~\sup_{\|x\| = 1}\| (\lambda T - I ) x \|^2 \\
\Rightarrow \inf _{\lambda \in C} \| (\lambda T - I ) y \|^2 & \leq & \inf _{\lambda \in C} ~~\sup_{\|x\| = 1}\| (\lambda T - I ) x \|^2
\end{eqnarray*}
This holds for all y in H with $ \|y\| = 1 $. So
\[ \sup_{\|x\| = 1}~~ \inf _{\lambda \in C} \| (\lambda T - I )x \|^2 \leq \inf _{\lambda \in C} ~~\sup_{\|x\| = 1}\| (\lambda T - I ) x \|^2 ~.\]
For the converse part we use the notion of center of mass of an operator. As T is a bounded linear operator  so there exists a  scalar $ \lambda_0  $ such that
\[ \|  \lambda_0 T - I \| = \inf _{\lambda \in C} \| (\lambda T - I ) \|.\]
Further there exists a sequence  $ \{x_{n}\},  \|x_{n}\| = 1 $ such that $ ((\lambda_0 T - I)x_{n}, Tx_{n})  \rightarrow  0 $ and $ \|(\lambda_0 T - I)x_{n}\| \rightarrow  \|\lambda_0 T - I\| $. \\
Now
\begin{eqnarray*}
\| \lambda_0 T - I \|^2 & = & \lim_{n\rightarrow \infty} \|(\lambda_0 T - I)x_n \|^2 \\
                         & = & \lim_{n\rightarrow \infty} \{ 1 - 2 Re \lambda_0    \langle T x_n, x_n \rangle + {\lambda_0}^2 \|Tx_n\|^2 \} \\
                         & = & \lim_{n\rightarrow \infty} \{ 1 - \frac{{ \mid \langle T x_n, x_n \rangle \mid}^2}{ \|Tx_n\|^2 } \}\\
                         & \leq & \sup_{\|x\| = 1} \{ 1 - \frac{{ \mid \langle T x, x \rangle \mid }^2}{ \|Tx\|^2 } \}~~~~~(*)
\end{eqnarray*}
For a fixed y in H with $ \|y \| = 1 $ we see that
\[ \| (\lambda T - I ) y \|^2  =  1 - 2 Re \lambda   \langle T y, y \rangle + {\lambda}^2 \|Ty\|^2 \]
 \[ \Rightarrow \| (\lambda T - I ) y \|^2  = \|Ty\|^2 \{ \mid \lambda - \frac{ \langle y,Ty \rangle}{\|Ty\|^2} \mid ^2 \} +   \{ 1 - \frac{ \mid \langle Ty,y \rangle \mid ^2}{ \|Ty\|^2 }\} .\]
achieves its minimum at $ \lambda(y) = \frac{ \langle y,Ty \rangle}{\|Ty\|^2} $ and the  minimum value is
\[ \inf_{\lambda \in C } \| (\lambda T - I ) y \|^2 = 1 - \frac{{ \mid \langle T y, y \rangle \mid}^2}{ \|Ty\|^2 } .\]
Using this in (*) we get
\[ \| \lambda_0 T - I \|^2 \leq \sup_{\|x\| = 1} \inf_{\lambda \in C} \| (\lambda T - I ) x \|^2 \]
and so
\[ \inf _{\lambda \in C} ~~\sup_{\|x\| = 1}\| (\lambda T - I ) x \|^2 \leq \sup_{\|x\| = 1}~~ \inf _{\lambda \in C} \| (\lambda T - I ) x \|^2 .\]
Thus we have the  equality
\[ \sup_{\|x\| = 1}~~ \inf _{\lambda \in C} \| (\lambda T - I ) x \|^2 = \inf _{\lambda \in C} ~~\sup_{\|x\| = 1}\| (\lambda T - I ) x \|^2 ~.\]
 \noindent \textbf{Theorem 3.7} \\
 \noindent Suppose $ \lambda_0 $ is the total center of mass of a  bounded linear operator T. Then, $\left|\cos\right| T = \lim_{n \rightarrow \infty} \frac{\mid (Tx_{n},x_{n})\mid}{\|Tx_{n}\|}$ where $\{ x_{n} \}$ is a sequence of unit vectors in H, $((I- \lambda_{0}T)x_{n},Tx_{n}) \rightarrow 0$ and $\|(I-\lambda_{0}T)x_{n}\| \rightarrow \|I-\lambda_{0}T\|$\\
 \noindent \emph{\bf Proof.}
 In \cite{18} we proved that if  $ \lambda_0  $ is the total center of mass of T then there exists a sequence  $ \{x_{n}\},  \|x_{n}\| = 1 $ such that $ ((I - \lambda_0 T )x_{n}, Tx_{n})  \rightarrow  0 $ and $ \|(I - \lambda_0 T )x_{n}\| \rightarrow  \|I - \lambda_0 T\| $ and
\[ \|I - \lambda_0 T\| = \inf_{\lambda \in C} \| I - \lambda T \| = \inf _{\lambda \in C} ~~\sup_{\|x\| = 1}\| (I - \lambda T  ) x \| ~.\]
As in Theorem 3.6  we have
\begin{eqnarray*}
\| \lambda_0 T - I \|^2 & = & \lim_{n\rightarrow \infty} \{ 1 - \frac{{ \mid \langle T x_n, x_n \rangle \mid}^2}{ \|Tx_n\|^2 } \}\\
                         & \leq & \sup_{\|x\| = 1} \{ 1 - \frac{{ \mid \langle T x, x \rangle \mid}^2}{ \|Tx\|^2 } \} \\
                         & \leq & \sup_{\|x\| = 1} \inf_{\lambda \in C } \| (I - \lambda T  ) x \|^2 \\
                         & = & \| \lambda_0 T - I \|^2
\end{eqnarray*}
This shows that \[ \inf_{\|x\| = 1 }  \frac{{ \mid \langle T x, x \rangle \mid }}{ \|Tx\| } = \lim_{n\rightarrow \infty} \frac{{ \mid \langle T x_n, x_n \rangle \mid}}{ \|Tx_n\| } = \left|\cos \right| T\]
\\
We next give an example to calculate total antieigenvalue for a finite dimensional operator using the above method.\\
\noindent \textbf{Example 3.8} \\
$T=\left(
  \begin{array}{cc}
    1 & 0 \\
    0 & 1+i \\
  \end{array}
\right)
$ be an operator on a two dimensional complex Hilbert space H. Let, $z=(z_{1}, z_{2})^{t} \in H$, where $\left|z_{1}\right|^{2} + \left| z_{2}\right|^{2} = 1$. Then, $\left|(Tz,z)\right|=\sqrt{1+\left|z_{2}\right|^{4}}$ and $\|Tz\|=\sqrt{ 1+\left|z_{2}\right|^{2}}$. Now, $\sup_{\|z\| = 1}\inf_{\lambda \in C}\|(\lambda T - I)z\| =\sqrt{2}-1$ and this supremum is attained by the vector $z_{0}=(z_{1},z_{2})^{t}$, where $\left|z_{2}\right|^{2}=\sqrt{2}-1$. Then, $\lambda_{0}=\frac{1-i(\sqrt{2}-1)}{\sqrt{2}}$ and $\|\lambda_{0}T-I\|=\sqrt{2}-1$. Now, $\left|\cos\right| T = \frac{Re(Tz_{0},z_{0})}{\|Tz_{0}\|}=\sqrt{2\sqrt{2}-2}$.\\

\noindent \textbf{Acknowledgement}: We would like to thank Professor T. K. Mukherjee and Professor K. C. Das for their suggestion while preparing the paper. Second author would like to thank CSIR, India for supporting his research.

{\normalsize Kallol Paul}\\
{\normalsize Department of Mathematics}\\
{\normalsize Jadavpur University, Kolkata 700032,
West Bengal, INDIA}\\
{\normalsize email:} \emph{\normalsize kalloldada@yahoo.co.in}{\normalsize }\\
{\normalsize }\\

\noindent {\normalsize Gopal Das}\\
{\normalsize Department of Mathematics, Jadavpur University, Kolkata
700032, INDIA }\\
{\normalsize email:} \emph{\normalsize gopaldasju@gmail.com}{\normalsize }\\
 \end{document}